# Additive arithmetic functions with limit normal distribution

Victor Volfson

ABSTRACT. This paper proves several assertions on sufficient conditions for the convergence of additive arithmetic functions to the normal distribution.

A generalization of the Erdos-Kac theorem was proved and determines the rate of convergence of additive arithmetic functions to the normal distribution in the cases considered.

The paper provides examples of using the proved assertions.

## 1. INTRODUCTION

An arithmetic function is generally defined as a function $f(m), m=1,...,n$ on the set of natural numbers and taking values on the set of complex numbers. The name arithmetic function is related to the fact that this function expresses some arithmetic property of the natural series.

Let $m_1, m_2$ are relatively prime natural numbers.

An arithmetic function $f(m), m=1,...,n$ is additive if it satisfies:

$$f(m_1 m_2) = f(m_1) + f(m_2), f(1) = 0. \qquad (1.1)$$

Let $p$ is an arbitrary prime number and $\alpha > 1$ is a natural number and $m = \prod_{p^\alpha | m} p^\alpha$, then based on (1.1) the following holds for the additive arithmetic function $f(m), m=1,...,n$:

$$f(m) = \sum_{p^\alpha | m} f(p^\alpha). \qquad (1.2)$$

An additive arithmetic function $f(m), m=1,...,n$ is strongly additive if:

$$f(p^\alpha) = f(p). \qquad (1.3)$$

Based on (1.2), (1.3), the property is satisfied for strongly additive functions - $f(m) = \sum_{p|m} f(p)$.

The axiomatics of probability theory provides the possibility of using probability theory to find the distribution of any arithmetic functions. Any initial interval of the natural series $[1,n]$ can be naturally transformed into a probability space $(\Omega_n, A_n, P_n)$ by taking for $\Omega_n = \{1,2,...,n\}$, $A_n$ - all subsets $\Omega_n$, $P_n(A) = \#(m \in A)/n$, where $\#(m \in A)$ is the number of natural numbers in the subset $A$.

Then an arbitrary real arithmetic function $f(m), m=1,...,n$ on $\Omega_n$ can be considered as a random variable $x_n$ on this probability space:

$$x_n(m) = f(m)(1 \le m \le \text{n}). \qquad (1.4)$$

In particular, we can talk about the average value of an arithmetic function (mathematical expectation):

$$E[x_n] = E[f,n] = \frac{1}{n}\sum_{m=1}^{n} f(m)\,, \tag{1.5}$$

variance:

$$D[x_n] = D[f,n] = \frac{1}{n}\sum_{m=1}^{n} f^2(m) - \Big(\frac{1}{n}\sum_{m=1}^{n} f(m)\Big)^2\,, \tag{1.6}$$

and for a real arithmetic function $f(m), m = 1,...,n$ - about the distribution function:

$$F_n(y) = P_n\{m \le n : f(m) \le y\} \tag{1.7}$$

and the characteristic function:

$$\varphi_{x_n}(t) = \varphi_{f,n}(t) = E[\mathrm{e}^{\mathrm{i}tf_n}] = \frac{1}{n}\sum_{m=1}^{n} \mathrm{e}^{\mathrm{i}tf(m)}\,. \tag{1.8}$$

Erdős and Wintner proved [1] that the condition of convergence of series of a real additive arithmetic function $f(m), m = 1,...,n$ for a positive value $R$:

$$\sum_{|f(p)|\le R} \frac{f(p)}{p}\,,$$

$$\sum_{|f(p)|\le R} \frac{f^2(p)}{p}\,,$$

$$\sum_{|f(p)|> R} \frac{1}{p}\,.$$

is sufficient for the convergence of a sequence of distribution functions $\{F_n(y)\}$ to the limit distribution function $F(y)$ as $n \to \infty$. The idea of the proof based on theorems on sums of independent random variables.

It is interesting the proof of the fact of convergence of the distribution function to the limit distribution function and finding the asymptotic distribution laws for additive arithmetic functions in the case when the specified series diverge.

Erdős and Katz [2] used the central limit theorem to obtain that the limit distribution law for a normalized strongly additive arithmetic function $\omega(m), m = 1,...,n$ is normal. They later proved that this holds for all strongly additive arithmetic functions for which $B(n) = \sqrt{\sum_{p \le n} \frac{f^2(p)}{p}} \to \infty$ as n $\to \infty$ and $| f(p) | \le 1$.

Kubilius [3] introduced a class H of additive arithmetic functions. This class includes real additive arithmetic functions $f(m), m = 1,...,n$ for which $D(n) \to \infty$ as $n \to \infty$ and, there is an infinitely increasing function $r = r(n)$ such that $\ln r(n) / \ln(n) \to 0, D(r(n)) / D(n) \to 1, n \to \infty$, where $D(n) = \sqrt{\sum_{p^k \le n} \frac{f^2(p^k)}{p^k}}$. These conditions become $B(n) \to \infty$ as $n \to \infty$ and $\ln r(n) / \ln(n) \to 0, B(r(n)) / B(n) \to 1, n \to \infty$ in the case of strongly additive arithmetic functions.

It was proved that for strongly additive arithmetic functions of class H $\ln B(n) = o(\ln \ln n)$ is satisfied as $n \to \infty$.

In addition, the following assertion was proved for strongly additive arithmetic functions of class H.

Assertion 1

The distribution laws

$$P_n\{\frac{f(m) - A(n)}{B(n)} < x\},$$

where $f(m), m = 1,...,n$ is a real strongly additive arithmetic function from class H and $A(n) = \sum_{p \le n} \frac{f(p)}{p}$, converge to the limit with variance 1 if and only if there is a non-decreasing function $K(u)$ with variation 1 for which at all points of continuity $K(u)$ at $n \to \infty$ the following holds:

$$\frac{1}{B^2(n)} \sum_{\substack{p \le n \\ f(p) > uB(n)}} \frac{(f(p))^2}{p} \to K(u). \tag{1.9}$$

The logarithm of the characteristic function of the limit law is calculated using the Kolmogorov formula:

$$\ln\varphi(t)=\int_{-\infty}^{\infty}\frac{e^{itu}-1-itu}{u^2}dK(u).$$

The function $K(u)=0$ for $u\le 0$ and $K(u)=1$ for $u>0$ in the case of a normal distribution, therefore condition (1.9) goes over to (1.10) and in combination with Assertion 1 gives Assertion 2 (analogous to the Lindeberg condition).

Assertion 2

Let $f(m), m=1,...,n$ is a real strongly additive arithmetic function. If $B(n)\to\infty$ and for all fixed $\xi>0$ the condition is satisfied:

$$\frac{1}{B^2(n)}\sum_{\substack{p\le n\\ f(p)>\xi B(n)}}\frac{(f(p))^2}{p}\to 0, \tag{1.10}$$

then

$$P_n\{\frac{f(m)-A(n)}{B(n)}<x\}\to \mathrm{N}(0,1), n\to\infty, \tag{1.11}$$

where $A(n)=\sum_{p\le n}\frac{f(p)}{p}$.

The condition (1.10) is also necessary for the fulfillment of (1.11) in the case $f(m)\in H$.

Let us set the task of proposing simpler sufficient conditions for the fulfillment of (1.11) for additive arithmetic functions.

It is also interesting - the generalization of the above-mentioned Erdös–Kac theorem [4,5,6].

Another topical issue is determining the rate of convergence of additive arithmetic functions to a normal distribution. These issues were studied in [7,8,9]. Let me remind you that the standard rate of convergence is $O(1/\sqrt{\ln\ln n})$.

## 2. DETERMINATION OF MOMENTS OF ADDITIVE ARITHNETIC FUNCTIONS

We call additive arithmetic functions $f_r(m), m = 1,...,n$ "truncated" if they have only divisors $p^\alpha | m$ with $p \le r$.

Truncated additive arithmetic functions $f_r(m), m = 1,...,n$ can be represented as a sum of asymptotically independent random variables [3] at $n \to \infty$ on the probability space specified in Chapter 1:

$$f_r(m) = \sum_{p^\alpha \le r(n)} X_p + o(1), \tag{2.1}$$

where $X_p$ is equal to $f(p^\alpha)$ with probability $1/p^\alpha$ and is equal to 0 with the opposite probability.

Let us consider an arbitrary additive arithmetic function $f(m), m = 1,...,n$ ( $f(m) = \sum_{p^\alpha | m} f(p^\alpha)$ ).

Based on [3], the additive arithmetic function $f(m) \in H$ at $n \to \infty$ can be represented similarly to (2.1):

$$f(m) = \sum_{p^\alpha \le n} X_p + o(1). \tag{2.2}$$

Taking into account (2.2), the asymptotics of the moments $f(m) \in H$ coincides with the moments of the sum of independent random variables - $S(n) = \sum_{p^\alpha \le n} X_p$ .

This simplifies the application of probabilistic methods: it allows one to find moments of a random variable $S(n) = \sum_{p^\alpha \le n} X_p$ that correspond to the asymptotics of moments $f(m), m = 1,...,n$, to use the central limit theorem to find the limit distribution $f(m), m = 1,...,n$, etc.

Above we used formulas for the asymptotics of the mean value of a strongly additive arithmetic function $A(n)=\sum_{p\le n}\frac{f(p)}{p}$ and the asymptotics of the variance of a strongly additive arithmetic function $B^2(n)=\sum_{p\le n}\frac{f^2(p)}{p}$. Let us prove them.

The mathematical expectation $S(n)$ is calculated as the sum of the mathematical expectations of random variables $X_p$:

$$A(n)=E[\sum_{p^\alpha\le n}X_p]=\sum_{p^\alpha\le n}E[X_p]. \qquad (2.3)$$

Since $X_p$ takes the value $f(p^\alpha)$ with probability $1/p$ and 0 with probability $1-1/p$, then:

$$E[X_p]=f(p^\alpha)\cdot\frac{1}{p}+0\cdot(1-\frac{1}{p})=\frac{f(p^\alpha)}{p}. \qquad (2.4)$$

Substitute (2.4) into (2.3) and get:

$$A(n)=E[\sum_{p^\alpha\le n}X_p]=\sum_{p^\alpha\le n}E[X_p]=\sum_{p^\alpha\le n}\frac{f(p^\alpha)}{p}. \qquad (2.5)$$

Based on (2.5), we obtain for a strongly additive function:

$$A^*(n)=\sum_{p\le n}\frac{f(p)}{p},$$

which corresponds to the condition of Assertion 1.

The variance $S(n)$ is calculated as the sum of the variances $X_p$, since the random variables are independent:

$$D^2(n)=D^2[\sum_{p^\alpha\le n}X_p]=\sum_{p^\alpha\le n}D^2[X_p]. \qquad (2.6)$$

The variance $X_p$ is:

$$D^2[X_p]=E[{X_p}^2]-(E[X_p])^2. \qquad (2.7)$$

Since $X_p{}^2 = f(p^\alpha)^2$ with probability $\frac{1}{p}$ and 0 with probability $1-\frac{1}{p}$ then:

$$E[X_p{}^2] = f(p^\alpha)^2 \cdot \frac{1}{p} + 0 \cdot (1-\frac{1}{p}) = \frac{f(p^\alpha)^2}{p}. \tag{2.8}$$

Having in mind (2.4), (2.7) and (2.8):

$$D^2[X_p] = E[X_p{}^2] - (E[X_p])^2 = \frac{f(p^\alpha)^2}{p} - (\frac{f(p^\alpha)}{p})^2 = \frac{f(p^\alpha)^2}{p}(1-\frac{1}{p}). \tag{2.9}$$

Substituting (2.9) into (2.6) we obtain the variance $S(n)$:

$$D^2(n) = D^2[\sum_{p^\alpha \le n} X_p] = \sum_{p^\alpha \le n} D^2[X_p] = \sum_{p^\alpha \le n} \frac{f(p^\alpha)^2}{p}(1-\frac{1}{p}). \tag{2.10}$$

Expression (2.10) for large values $p$ can be simplified:

$$D^2(n) = \sum_{p^\alpha \le n} \frac{f^2(p^\alpha)}{p}. \tag{2.11}$$

Based on (2.11), we obtain for a strongly additive function:

$$B^2(n) = \sum_{p \le n} \frac{f^2(p)}{p},$$

which corresponds to the condition of Assertion 1.

Similarly to (2.11) we find the third absolute moment:

$$\rho(n) = E[|X_p|^3] = \sum_{p^\alpha \le n} \frac{|f(p^\alpha)|^3}{p}. \tag{2.12}$$

Assertion 3

Let there is an additive arithmetic function $f(m) \in H$ and a truncated additive arithmetic function $f_r(m) \in H$. Let us denote the asymptotics of the mean value $f_r(m)$ by - $A_r(n)$ and the asymptotics of the variance $f_r(m)$ by - $D_r{}^2(n)$. Then, for $n \to \infty$, the following holds:

$$A(n) \sim A_r(n) \text{ and } D^2(n) \sim D_r{}^2(n) \tag{2.13}$$

Proof

Having in mind the definition of class H of additive arithmetic functions given in Chapter 1 at $n \to \infty$:

$$D^2(n) \sim D_r^{\,2}(n).$$

Based on Lemmas 4.1 and 4.2 of [3], for $n \to \infty$:

$$A(n) \sim A_r(n).$$

Consequence 4

Considering that a strongly additive arithmetic function $f^*(m)$ is a special case of a truncated additive arithmetic function $f_r(m)$, then it follows from Assertion 2:

$$D^2(n) \sim B^2(n) \text{ and } A(n) \sim A^*(n) \text{ for } n \to \infty. \tag{2.14}$$

Since functions of class H can be represented as a sum of independent random variables, this allows one to apply classical results of probability theory, such as the Berry-Essen inequality to estimate the rate of convergence to a normal distribution.

The rate of convergence in the Berry–Essen inequality for the sum of independent random variables with different distributions $S(n)$ is estimated as:

$$\Delta = \sup | P(\frac{S(n)}{D(n)} \le x) - \Phi(x) | = O(\frac{\rho(n)}{D^3(n)}), \tag{2.15}$$

where $\Phi(x)$ is the standard normal distribution function.

## 3. SUFFICIENT CONDITIONS FOR CONVIRGENCE TO THE NORMAL DISTRIBUTION FOR ADDITIVE ARITHMETIC FUNCTIONS

Let us consider simpler and more convenient formulations of sufficient conditions for the fulfillment of (1.11) for additive arithmetic functions.

Assertion 5

Let $f(m), m=1,...,n$ is a real strongly additive arithmetic function, i.e. $f(m)=\sum_{p|m} f(p)$ and the condition is satisfied at $n \to \infty$ in this case:

$$\max_{p\le n} |f(p)| = o(B(n)) \tag{3.1}$$

Then

$$P_n\{\frac{f(m)-A(n)}{D(n)} < x\} \to \Phi(x), n \to \infty \tag{3.2}$$

Proof

Based on Assertion 2 (see introduction), the condition must be satisfied for all fixed $\xi > 0$:

$$\frac{1}{B^2(n)} \sum_{\substack{p\le n \\ f(p)>\xi B(n)}} \frac{(f(p))^2}{p} \to 0. \tag{3.3}$$

Condition (3.3) can be written as $n \to \infty$:

$$\frac{1}{B^2(n)} \sum_{\substack{p\le n \\ f(p)>\xi B(n)}} \frac{(f(p))^2}{p} = \frac{1}{B^2(n)} \sum_{p\le n} \frac{(f(p))^2}{p} I(\frac{|f(p)|}{B(n)} > \varepsilon) \to 0, \tag{3.4}$$

where I(.) is the indicator function.

Since $B^2(n) = \sum_{p\le n} \frac{(f(p))^2}{p}$, then condition (3.4) is satisfied under condition (3.1).

Therefore, based on Assertion 2, the following holds:

$$P_n\{\frac{f(m)-A(n)}{D(n)} < x\} \to \Phi(x), n \to \infty. \tag{3.5}$$

Assertion 6

Let $f(m), m=1,...,n$ a is a real additive arithmetic function, i.e. $f(m) = \sum_{p^\alpha|m} f(p^\alpha)$. It belongs to the class H and the condition is satisfied at $n \to \infty$:

$$\max_{p\le n} | f(p) | = o(B(n)) \qquad (3.6)$$

Then $P_n\{\frac{f(m)-A(n)}{D(n)} < x\} \to \Phi(x), n\to\infty$ with the rate of convergence $\Delta = O(\frac{\rho(n)}{D^3(n)})$.

Proof

Since an additive arithmetic function $f(m)\in H$, then based on [3] there is such a real strongly additive arithmetic function $f^*(m)=\sum_{p|m} f(p)\in H$ for which the following holds:

$$P_n\{\frac{f(m)-A(n)}{D(n)} < x\} = P_n\{\frac{f(m)-A^*(n)}{B(n)} < x\}, n\to\infty, \qquad (3.7)$$

i.e. the limit distributions of the normalized values of these arithmetic functions coincide.

Therefore, all conditions of Assertion 5 are satisfied for $f^*(m)=\sum_{p|m} f(p)\in H$ and based on (3.2):

$$P_n\{\frac{f(^*m)-A(n)}{B(n)} < x\} \to \Phi(x), n\to\infty. \qquad (3.8)$$

Having in mind (3.7) and (3.8) the following holds:

$$P_n\{\frac{f(m)-A(n)}{D(n)} < x\} \to \Phi(x), n\to\infty. \qquad (3.9)$$

Based on (2.15) we get:

$$\Delta = \sup | P(\frac{S(n)}{D(n)} \le x) - \Phi(x) | = O(\frac{\rho(n)}{D^3(n)}) \text{ при } n\to\infty.$$

Therefore $P_n\{\frac{f(m)-A(n)}{D(n)} < x\} \to \Phi(x), n\to\infty$, with the rate of convergence $\Delta = O(\frac{\rho(n)}{D^3(n)})$.

As an example of using Assertion 6, we consider the additive arithmetic function $f(m)=\ln\tau_k(m)$, where $\tau_k(m)$ is the number of ways to represent $m$ as a product of $k$ factors (taking into account the order).

It is performed for prime $p$ and its degree $\alpha$: $\ln(\tau_k(p^\alpha)) = \ln(C^{k-1}{}_{a+k-1})$.

Let us check whether the additive arithmetic function $f(m) = \ln\tau_k(m)$ belongs to class H.

The variance of an arithmetic function $f(m) = \ln\tau_k(m)$ tends to:

$$D^2(n) \sim \ln^2(k)\ln\ln n \to \infty$$

at $n \to \infty$.

Let's choose $r(n) = \exp\dfrac{\ln n}{(\ln\ln n)^a}, a > 0$ and check the conditions:

$$\frac{\ln r(n)}{\ln n} = \frac{1}{(\ln\ln n)^a} \to 0 \text{ when } n \to \infty.$$

$$\ln\ln r(n) \sim \ln\ln n,\ D(r(n)) \sim \ln k\sqrt{\ln\ln r(n)} \sim \ln k\sqrt{\ln\ln n}.$$

$$\frac{D(r(n))}{D(n)} \sim \frac{\ln k\sqrt{\ln\ln n}}{\ln k\sqrt{\ln\ln n}} = 1.$$

Thus, all conditions are met $f(m) = \ln\tau_k(m) \in H$ and $f(m) = \ln\tau_k(m) \in H$.

Therefore, the formulas of Section 2 can be used.

Based on (2.5), asymptotics of the mean value $f(m) = \ln\tau_k(m)$:

$$A(n) = \sum_{p^\alpha \le n} \frac{f(p^\alpha)}{p} \sim \ln(k)\ln\ln n.$$

Based on (2.10), the asymptotics of the variance of the arithmetic function $f(m) = \ln\tau_k(m)$:

$$D^2(n) = \sum_{p^\alpha \le n} \frac{f(p^\alpha)^2}{p}\left(1 - \frac{1}{p}\right) \sim \ln^2(k)\ln\ln n.$$

Based on (2.12), the asymptotics of the third absolute moment $f(m) = \ln\tau_k(m)$:

$$\rho(n) = \sum_{p^\alpha \le n} \frac{|f(p^\alpha)|^3}{p} \sim \ln^3(k)\ln\ln n.$$

Now let us check the fulfillment of the sufficient condition (3.6):

$$\max_{p\le n} |f(p)| = \max_{p\le n} |\ln \tau_k(p)| = \max_{p\le n} |\ln C_k^{k-1}| = \ln k\,,$$

$$B(n) = \sqrt{\sum_{p\le n} \frac{f(p)^2}{p}(1-\frac{1}{p})} \sim \ln(k)\sqrt{\ln\ln n}\,,$$

$$\max_{p\le n} |f(p)| = \ln k = o(B(n)) = o(\ln k\sqrt{\ln\ln n})\,.$$

Since the sufficient condition (3.6) is satisfied, then

$$P_n\{\frac{f(m)-A(n)}{D(n)} < x\} = P_n\{\frac{\ln(\tau_k(m)) - \ln k\ln\ln(n)}{\ln k\sqrt{\ln\ln(n)}} < x\} \to \Phi(x), n \to \infty\,.$$

Let us determine the rate of convergence at $n \to \infty$ for this case:

$$\Delta = O(\frac{\rho(n)}{D^3(n)}) = O(\frac{\ln^3(k)\ln\ln n}{\ln^3(k)(\ln\ln n)^{3/2}}) = O(\frac{1}{\sqrt{\ln\ln n}})\,,$$

i.e. it coincides with the standard one.

Let us find less stringent sufficient conditions for the fulfillment of (1.11) for additive arithmetic classes H than the condition in Assertion 6.

Assertion 7

Let $f(m), m = 1,...,n$ is a real additive arithmetic function, i.e. $f(m) = \sum_{p^\alpha | m} f(p^\alpha)$. It belongs to the class H and the condition is satisfied at $n \to \infty$:

$$\max_{p\le n} |f(p)| \le \varepsilon_n B(n)), \varepsilon_n \to 0\,. \tag{3.10}$$

Тогда $P_n\{\frac{f(m)-A(n)}{D(n)} < x\} \to \Phi(x), n \to \infty$ со скоростью сходимости $\Delta = O(\frac{\rho(n)}{D^3(n)})$.

Then $P_n\{\frac{f(m)-A(n)}{D(n)} < x\} \to \Phi(x), n \to \infty$ with the rate of convergence $\Delta = O(\frac{\rho(n)}{D^3(n)})$.

Proof

Let us divide both sides of condition (3.10) by $B(n)$:

$$\frac{\max_{p\le n}|f(p)|}{B(n)} \le \varepsilon_n .$$

Since by the condition $\varepsilon_n \to 0$ at $n \to \infty$, then passing to the limit we obtain:

$$\lim_{n\to\infty} \frac{\max_{p\le n}|f(p)|}{B(n)} \le \lim_{n\to\infty} \varepsilon_n = 0 .$$

Thus:

$$\lim_{n\to\infty} \frac{\max_{p\le n}|f(p)|}{B(n)} = 0 .$$

Based on the definition of o-small:

$\max_{p\le n}|f(p)| = o(B(n))$ when $n \to \infty$.

Therefore, sufficient condition (3.6) of Assertion 6 is satisfied.

Therefore, based on Assertion 6: $P_n\{\frac{f(m)-A(n)}{D(n)} < x\} \to \Phi(x), n \to \infty$ with the rate of convergence $\Delta = O(\frac{\rho(n)}{D^3(n)})$

The sufficient condition $\max_{p\le n}|f(p)| \le \varepsilon_n B(n)), \varepsilon_n \to 0$ is less stringent than $\max_{p\le n}|f(p)| = o(B(n))$ for $n \to \infty$, but still guarantees normal convergence. It allows individual contributions $|f(p)|$ to grow, but at a controlled rate.

Let's look at an example for Assertion 7.

Suppose $f(p) = \ln\ln p$. It can be verified that the strongly additive function $f(m) = \sum_{p|m} \ln\ln p \in H$. We will prove this later in a more general case (see Assertion 9).

Now let's take $\varepsilon_n = \sqrt{\frac{3}{\ln\ln n}}$, then condition (3.10) will be written as:

$$\max_{p\le n}|f(p)| \le \frac{\sqrt{3}B(n)}{\sqrt{\ln\ln n}} . \qquad (3.11)$$

In this case

$$B(n) = \sqrt{\sum_{p \le n} \frac{f(p)^2}{p}(1 - \frac{1}{p})} = \sqrt{\sum_{p \le n} \frac{(\ln\ln p)^2}{p}(1 - \frac{1}{p})} \sim \frac{(\ln\ln n)^{3/2}}{\sqrt{3}}.$$

Поэтому

$$\max_{p \le n} | f(p) | = \ln\ln n \le \frac{(\ln\ln n)^{3/2}\sqrt{3}}{\sqrt{3}(\ln\ln n)^{1/2}} = \frac{\sqrt{3}B(n)}{\sqrt{\ln\ln n}},$$

что соответствует (3.11).

Therefore,

$$\max_{p \le n} | f(p) | = \ln\ln n \le \frac{(\ln\ln n)^{3/2}\sqrt{3}}{\sqrt{3}(\ln\ln n)^{1/2}} = \frac{\sqrt{3}B(n)}{\sqrt{\ln\ln n}},$$

which corresponds to (3.11).

Thus, Assertion 7 is satisfied, therefore:

$$P_n\{\frac{f(m) - A(n)}{D(n)} < x\} = P_n\{\frac{\sum_{p|m} \ln\ln p - \frac{(\ln\ln n)^2}{2}}{\frac{(\ln\ln n)^{3/2}}{\sqrt{3}}} < x\} \to \Phi(x), n \to \infty.$$

Next, we will consider a more general case and the rate of convergence to the normal distribution determine for it (see Assertion 9).

Now let us assume that condition (3.10) is satisfied almost everywhere, i.e. exceptions to the fulfillment of this condition are possible. Let us ask ourselves how many exceptions are possible in the case, in which normal convergence is still fulfilled.

Assertion 8

Let $f(m), m = 1,...,n$ is a real additive arithmetic function, i.e. $f(m) = \sum_{p^\alpha | m} f(p^\alpha)$, it belongs to class H and two conditions are satisfied at $n \to \infty$:

1. A constraint on the maximum contribution, similar to (3.10)

$\max_{p\le n} | f(p) |\le \varepsilon_n B(n)), \varepsilon_n \to 0$.

2. Control over the number of exceptions

$$\#\{p \le n : | f(p) |> \varepsilon_n B(n)\} = o(1/\varepsilon^2{}_n). \qquad (3.12)$$

Then $P_n\{\frac{f(m) - A(n)}{D(n)} < x\} \to \Phi(x), n \to \infty$ with the rate of convergence $\Delta = O(\frac{\rho(n)}{D^3(n)})$.

Proof

Let's divide prime numbers $p \le n$ into two groups:

$| f(p) |\le \varepsilon_n B(n)$,

$| f(p) |> \varepsilon_n B(n)$.

By condition (3.12), the number of primes in the second group is:

$$\#\{p \le n : | f(p) |> \varepsilon_n B(n)\} = o(1/\varepsilon^2{}_n). \qquad (3.13)$$

According to condition (3.10), the following is satisfied for each $p$ of the second group:

$| f(p) |> \varepsilon_n B(n)$.

Their total contribution to the variance of the arithmetic function:

$$\sum_{\substack{p\le n| \\ f(p)|>\varepsilon_n B(n)}} \frac{f^2(p)}{p} \le \frac{(\varepsilon_n B(n))^2}{\min p} \#\{p \le n : | f(p) |> \varepsilon_n B(n)\}. \qquad (3.14)$$

Having in mind (3.9), (3.14) and $p \ge 2$ we obtain:

$$\sum_{\substack{p\le n \\ |f(p)|>\varepsilon_n B(n)}} \frac{f^2(p)}{p} = o(\frac{1}{\varepsilon_n{}^2})(\varepsilon_n B(n))^2 = o(B^2(n)). \qquad (3.15)$$

Based on (3.15), taking into account that $f(m), m = 1,...,n$ belongs to class H, $B(n) \to \infty$ for $n \to \infty$ and $\varepsilon > \varepsilon_n$ we obtain:

$$\frac{1}{B^2(n)} \sum_{\substack{p \le n \\ f(p) > \xi B(n)}} \frac{(f(p))^2}{p} \to 0$$

Therefore, based on Assertion 2, $P_n\{\frac{f(m) - A(n)}{D(n)} < x\} \to \Phi(x), n \to \infty$

Having in mind (2.15) we obtain an estimate for the rate of convergence:

$$\Delta = \sup |P(\frac{S(n)}{D(n)} \le x) - \Phi(x)| = O(\frac{\rho(n)}{D^3(n)}) \text{ at } n \to \infty .$$

Let's look at an example for Assertion 8.

There is no normal convergence for a strongly additive arithmetic function $f(m) = \sum_{p|m} \ln\ln p \in H$ when introducing exceptions.

Therefore, as an example, let us take a strongly additive arithmetic function, which for most primes $p \le n$ is equal to $f(p) = 1$, and for specially selected primes $p_k \le n$ is equal to $f(p_k) = (\ln\ln n)^{1/4}$ with the number of exceptions $K = \ln\ln\ln n$ .

It can be shown that this strongly additive arithmetic function belongs to the class H.

Let us determine the variance of this strongly additive arithmetic function using the formula:

$$B^2(n) = \sum_{p \le n} \frac{f^2(p)}{p} .$$

Let's divide this sum into two parts: the contribution of ordinary primes and the contribution of exceptions.

Contribution to the variance of ordinary primes:

$$\sum_{p \le n} \frac{1^2}{p} \sim \ln\ln n . \qquad (3.16)$$

Contribution of exceptions to variance:

$$\sum_{k=2}^{\ln\ln\ln n} \frac{((\ln\ln n)^{1/4})^2}{p_k} = \sum_{k=2}^{\ln\ln\ln n} \frac{(\ln\ln n)^{1/2}}{p_k}.$$

Having in mind, that $p_k \sim k\ln k$ we get:

$$\sum_{k=2}^{\ln\ln\ln n} \frac{(\ln\ln n)^{1/2}}{k\ln k} \le (\ln\ln n)^{1/2} \ln\ln\ln\ln\ln(n).$$

This contribution is negligible compared to (3.16), hence:

$$B(n) \sim \sqrt{\ln\ln n}. \qquad (3.17)$$

Let us assume $\varepsilon_n = \frac{1}{(\ln\ln n)^{1/4}}$. Then the constraint on the maximum contribution is satisfied:

$$\max_{p\le n} |f(p)| = (\ln\ln n)^{1/4} \le \frac{\sqrt{\ln\ln n}}{(\ln\ln n)^{1/4}}.$$

Now let's check the execution of exceptions:

$$\#\{p \le n : |f(p)| > \varepsilon_n B(n)\} = \#\{p \le n : |f(p)| > 1\} = \ln\ln\ln\ln n = o(1/\varepsilon^2{}_n) = o(\sqrt{\ln\ln n}).$$

Thus, all the conditions of Assertion 8 are satisfied, therefore it is satisfied for this strongly additive arithmetic function $P_n\{\frac{f(m)-A(n)}{D(n)} < x\} \to \Phi(x), n\to\infty$ with the rate of convergence $\Delta = O(\frac{\rho(n)}{D^3(n)})$.

In this case: $D(n) = B(n) = (\ln\ln n)^{1/2}$, $D^3(n) = B^3(n) = (\ln\ln n)^{3/2}$.

Let us define $A(n)$ for this strongly additive arithmetic function by the formula:

$$A(n) = \sum_{p\le n} \frac{f(p)}{p}.$$

Contribution of ordinary common shares to $A(n)$:

$$\sum_{p\le n}\frac{1}{p}\sim \ln\ln n\,. \tag{3.18}$$

Contribution of exceptions to $A(n)$:

$$\sum_{k=2}^{\ln\ln\ln n}\frac{(\ln\ln n)^{1/4}}{p_k}\le (\ln\ln n)^{1/4}\ln\ln\ln\ln\ln n\,.$$

This contribution is negligible compared to (3.18), so:

$$A(n)\sim \ln\ln n\,. \tag{3.19}$$

Let us define $\rho(n)$ for this strongly additive arithmetic function by the formula:

$$\rho(n)=\sum_{p\le n}\frac{f^3(p)}{p}\,.$$

Contribution of ordinary common shares to $\rho(n)$:

$$\sum_{p\le n}\frac{1}{p}\sim \ln\ln n\,. \tag{3.20}$$

Contribution of exceptions to $\rho(n)$:

$$\sum_{k=2}^{\ln\ln\ln n}\frac{(\ln\ln n)^{3/4}}{p_k}\le (\ln\ln n)^{3/4}\ln\ln\ln\ln\ln n\,.$$

This contribution is negligible compared to (3.20), so:

$$\rho(n)\sim \ln\ln n\,. \tag{3.21}$$

Using (3.17) and (3.21) we determine the rate of convergence:

$$\Delta = O\left(\frac{\rho(n)}{D^3(n)}\right)=O\left(\frac{\ln\ln n}{(\ln\ln n)^{3/2}}\right)=O\left(\frac{1}{\sqrt{\ln\ln n}}\right). \tag{3.22}$$

Based on (3.22), the standard rate of convergence is satisfied in this case.

## 4. GENERALIZATION OF THE ERDOS-KAC THEOREM

Assertion 9

It is fulfilled for a real strongly additive arithmetic function $f(m)=\sum_{p|m}(\ln\ln p)^{\alpha}, \alpha>0$:

$$P_n\{\frac{f(m)-A(n)}{B(n)}<x\}\to\Phi(x), n\to\infty,$$

where $A(n)\sim\frac{(\ln\ln(n))^{a+1}}{a+1}$, $B(n)\sim\frac{(\ln\ln(n))^{a+1/2}}{\sqrt{2a+1}}$ with the rate of convergence $\Delta=O(\frac{1}{\sqrt{\ln\ln n}})$.

Proof

First, let us prove that $f(m)=\sum_{p|m}(\ln\ln p)^{\alpha}\in H$.

Let's find the variance $f(m)=\sum_{p|m}(\ln\ln p)^{\alpha}$:

$$Var(n)\sim\int_{2}^{n}\frac{(\ln\ln t)^{2\alpha}dt}{t\ln t}\sim\frac{(\ln\ln n)^{2\alpha+1}}{2\alpha+1}.$$

Therefore, the asymptotics $B(n)$ is equal to:

$$B(n)=\sqrt{Var(n)}\sim\frac{(\ln\ln n)^{\alpha+0,5}}{\sqrt{2\alpha+1}}. \tag{4.1}$$

Based on (4.1): $B(n)\to\infty, n\to\infty$. The first condition is satisfied.

Let's check whether the second condition is met.

Let's choose $r(n)=n^{1/\ln\ln n}$.

$$\frac{\ln r(n)}{\ln n}=\frac{\ln n^{1/\ln\ln n}}{\ln n}=\frac{1}{\ln\ln n}\to 0 \text{ for } n\to\infty.$$

$$\ln\ln r(n)\sim\ln\frac{\ln n}{\ln\ln n}=\ln\ln n-\ln\ln\ln n,\ B(r(n))\sim\frac{(\ln\ln r(n))^{\alpha+0,5}}{\sqrt{2\alpha+1}}\sim\frac{(\ln\ln n)^{\alpha+0,5}}{\sqrt{2\alpha+1}}.$$

$$\frac{B(r(n))}{B(n)} \sim \frac{(\ln\ln n)^{\alpha+0,5}}{(\ln\ln n)^{\alpha+0,5}} = 1.$$

Thus, all conditions are met and $f(m) = \sum_{p|m} (\ln\ln p)^{\alpha} \in H$.

Let's check the fulfillment of the second condition of Assertion 5:

$$\max_{p\le n} |f(p)| = o(B(n)).$$

Having in mind (4.1), the following is indeed true:

$$\max_{p\le n} |(\ln\ln(p))^{\alpha}| = (\ln\ln n)^{\alpha} = o(B(n)).$$

Therefore, all the conditions of Assertion 6 are fulfilled in this case, therefore $P_n\{\frac{f(m)-A(n)}{B(n)} < x\} \to \Phi(x), n \to \infty$ with the rate of convergence $\Delta = O(\frac{\rho(n)}{B^3(n)})$.

Based on (2.5), the following holds for a strongly additive arithmetic function $f(m)$:

$$A(n) = \sum_{p\le n} \frac{f(p)}{p} = \sum_{p\le n} \frac{(\ln\ln p)^{\alpha}}{p} \sim \int_{t\le n} \frac{(\ln\ln t)^{\alpha} dt}{t\ln t} \sim \frac{(\ln\ln n)^{\alpha+1}}{\alpha+1}.$$

Having in mind (2.12) it is true for a strongly additive arithmetic function $f(m)$:

$$\rho(n) = \sum_{p\le n} \frac{|f(p)|^3}{p} = \sum_{p\le n} \frac{(\ln\ln(p))^{3\alpha}}{p} \sim \int_{t\le n} \frac{(\ln\ln t)^{3\alpha} dt}{t\ln t} \sim \frac{(\ln\ln n)^{3\alpha+1}}{3\alpha+1}. \quad (4.2)$$

Based on (4.1) and (4.2), we determine the rate of convergence in this case:

$$\Delta = O(\frac{\rho(n)}{B^3(n)}) = O(\frac{(\ln\ln n)^{3\alpha+1}}{(\ln\ln n)^{3\alpha+1,5}}) = O(\frac{1}{\sqrt{\ln\ln n}}). \quad (4.3)$$

Thus, the standard rate of convergence to the normal distribution also holds in this case.

Note that when $\alpha = 0$ the value is $f(p) = (\ln\ln p)^{\alpha} = 1$. Based on the Erdös-Kac theorem $P_n\{\frac{f(m)-A(n)}{B(n)} < x\} \to \Phi(x), n \to \infty$ and the rate of convergence is the same $\Delta = O(\frac{1}{\sqrt{\ln\ln n}})$.

Therefore, Assertion 9 is satisfied for $f(m) = \sum_{p|m} (\ln\ln p)^{\alpha}, \alpha \ge 0$.

Corollary 10

Let there is a real additive arithmetic function $f(m) \in H$ for which holds $f(p) = O((\ln\ln p)^{\alpha}), \alpha \geq 0$ for the corresponding strongly additive arithmetic function. Then holds $P_n\{\frac{f(m) - A^*(n)}{B(n)} < x\} \to \Phi(x)$ for $n \to \infty$, where $A^*(n) = \sum_{p \leq n} \frac{f(p)}{p}$, $B(n) = \sqrt{\sum_{p \leq n} \frac{f^2(p)}{p}}$ with the rate of convergence $\Delta = O(\frac{1}{\sqrt{\ln\ln n}})$.

Proof

Based on Lemmas 4.1 and 4.2 of [3], if $f(m) \in H$, then the distribution laws of normalized functions of an additive arithmetic function $f(m)$ and the corresponding strongly additive arithmetic function $f^*(m)$ coincide:

$$P_n\{\frac{f(m) - A(n)}{D(n)} < x\} = P_n\{\frac{f^*(m) - A^*(n)}{B(n)} < x\}, n \to \infty.$$

Taking this into account, based on Assertion 9, the following holds: $P_n\{\frac{f(m) - A^*(n)}{B(n)} < x\} \to \Phi(x)$ for $n \to \infty$, where $A^*(n) = \sum_{p \leq n} \frac{f(p)}{p}$, $B(n) = \sqrt{\sum_{p \leq n} \frac{f^2(p)}{p}}$ with the rate of convergence $\Delta = O(\frac{1}{\sqrt{\ln\ln n}})$.

Now we consider a case that does not satisfy Corollary 10.

Let there is a strongly additive arithmetic function $f(m) = \sum_{p|m} f(p)$, where $f(p) = \ln^a p, a > 0$.

Let's find the variance $f(m) = \sum_{p|m} f(p)$:

$$B^2(n) \sim \int_{t \leq n} \frac{\ln^{2a} t dt}{t \ln t} \sim \frac{\ln^{2a} n}{2a}.$$

Let's check the condition $\ln B(n) = o(\ln\ln n)$:

$$\ln B(n) \sim a \ln\ln n - 0,5\ln(2a) = O(\ln\ln n).$$

Thus, this condition is not satisfied, therefore this strongly additive function does not belong to the class H. Consequently, $f(m) = \sum_{p|m} f(p)$, where $f(p) = \ln^a p, a > 0$ does not satisfy Corollary 10.

## 5. CONCLUSION AND SUGGESTIONS FOR FURTHER WORK

Next article will continue to study the asymptotic behavior of some arithmetic functions.

## 6. ACKNOWLEDGEMENTS

Thanks to everyone who has contributed to the discussion of this paper. I am grateful to Milan Pasteka of the Mathematical Institute, member of the Slovak Academy of Sciences, who read the full text of the article before its publication in the Archive.